\documentstyle[12pt]{article}
\titlepage

\title{On Parallel Sections of a Vector Bundle}
\author{  Richard Atkins \\
        richard.atkins@twu.ca \\
		Department of Mathematics\\ 
		Trinity Western University \\
		7600 Glover Road \\
		Langley, BC, V2Y 1Y1 Canada}
\date{}
\newtheorem{fact}{Fact}

\newtheorem{theorem}[fact]{Theorem}
\newtheorem{corollary}[fact]{Corollary}
\begin{document}
\maketitle
\begin{abstract}
We consider when a smooth vector bundle endowed with a connection possesses non-trivial, local parallel sections. This is accomplished by means of a derived flag of subsets of the bundle. The procedure is algebraic and rests upon the Frobenius Theorem.
\end{abstract}

\newpage
\section{Introduction}
A connection on a vector bundle is a type of differentiation that acts on vector fields.
Its importance lies in the fact that given a piecewise continuous curve connecting two 
points on the underlying manifold, the connection defines a linear isomorphism between the respective 
fibres over these points.  A renowned theorem of differential geometry states that when the Riemann curvature tensor of the connection vanishes, there exist local frames comprised of parallel sections. This paper presents a refinement of this result. That is, given a connection on a vector bundle we determine when there exist local parallel sections and we find the  subbundle they generate. This is accomplished by means of an algebraic construction of a derived flag of subsets
of the original vector bundle.

This question has been considered in the real analytic case by Trencevski (cf. \cite{dd}). 
Our solution, by contrast, follows from the Frobenius Theorem, which applies to smooth ($C^{\infty}$) data. 
Furthermore, while Trencevski's method relies upon power series expansions we take a more geometric approach to the problem.

By applying our methods to  the vector bundle of symmetric two-tensors over a manifold we obtain a solution to the problem
of determining when a connection is locally a metric connection. For the case of surfaces, this has also been dealt with in \cite{aa}.

\section{The Existence of Parallel Sections}

Let $\pi: W\rightarrow M$ be a smooth vector bundle and $W'$ a subset of $W$ with the 
following two properties: \\
P1: For each $x\in M$,  $W_{x}\cap W'$ is a linear subspace of $W_{x}:=\pi^{-1}(x)$. \\
P2: For each $ w\in W'$ there exists an open neighbourhood $U$ of $\pi(w)$ in $M$
and a smooth local section $X: U\subseteq M\rightarrow W'\subseteq W$ such that
$w = X(\pi(w))$. 

Let \[ \nabla :{\cal A}^{0}(W) \rightarrow {\cal A}^{1}(W) \]
be a connection on $W$, where ${\cal A}^{n}(W)$ denotes the space of local sections 
$U\subseteq M \rightarrow W \otimes \Lambda^{n}M$. Define a map 
\[ \widetilde{\alpha}: {\cal A}^{0}(W) \rightarrow {\cal A}^{1}(W/W') \] 
by \[ \widetilde{\alpha} : = \phi \circ \nabla, \]
where $W/W'$ is the quotient of $W$ and $W'$ taken fibrewise and
$\phi:W\otimes T^{*}M \rightarrow (W/W')\otimes T^{*}M$ denotes the natural projection.
For any local section $X:U\subseteq M \rightarrow W'\subseteq W$ and 
differentiable function $f:U\rightarrow \Re$
we have $\widetilde{\alpha}(fX) = f\widetilde{\alpha}(X)$. Thus, there corresponds to 
$\widetilde{\alpha}$ a map 
\[ \alpha_{W'} : W' \rightarrow (W/W')\otimes T^{*}M \]
acting linearly on each fibre of $W'$.
$\alpha_{W'}$ is the {\it second fundamental 1-form} of $W'$. 

Let $V$ be any subset of $W$ satisfying P1. Define ${\cal S}(V)$ to be the subset of
$V$ consisting of all elements $v$ for which there exists a smooth local section
$X:U\subseteq M \rightarrow V\subseteq W$ such that $v = X(\pi(v))$. Then 
${\cal S}(V)$ satisfies both P1 and P2.

We seek to construct the maximal flat subset $\widetilde{W}$,  of $W$. 
$\widetilde{W}$ may be obtained as follows. Set
\[ \begin{array}{lll}
   V^{(0)} & := & \{ w \in W \hspace{.03in} | \hspace{.03in} R(,)(w) = 0 \} \\
   W^{(i)}   & := & {\cal S}(V^{(i)}) \\
   V^{(i+1)} & := & ker \hspace{.03in} \alpha_{W^{(i)}} 
  \end{array} \]
where $R:TM\otimes TM\otimes W\rightarrow W$ 
denotes the curvature tensor of $\nabla$. This gives a sequence
\[ W \supseteq W^{(0)} \supseteq W^{(1)} \supseteq \cdots \supseteq W^{(k)} 
\supseteq \cdots \]
of subsets of $W$. Note that $W^{(i)}$ is not necessarily a vector bundle over $M$ since
the dimension of the fibres may vary from point to point.
For some $k \in { N}$, $W^{(l)}= W^{(k)}$ for all $l \geq k$. 
Define $\widetilde{W} = W^{(k)}$, with projection $\tilde{\pi}:\widetilde{W} \rightarrow M$. 

In order to extract information from $\widetilde{W}$ we need some concept of regularity. 
Accordingly, we say that the connection $\nabla$ is {\it regular at $x\in M$} if there exists 
a neighbourhood $U$ of $x$ such that $\tilde{\pi}^{-1}(U) \subseteq \widetilde{W}$
is a vector bundle over $U$. $\nabla$ is $regular$ if $\widetilde{W}$ is a vector bundle over
$M$. The dimension of the fibres of $\widetilde{W}$, for regular $\nabla$,
shall be denoted $rank \hspace{0.03in} \widetilde{W}$. 

\begin{theorem} \label{theorem:flag}
Let $\nabla$ be a connection on the smooth vector bundle $\pi:W\rightarrow M$. \\
(i) If $X:U \subseteq M \rightarrow W$ is a local parallel section then the image of $X$ lies in
$\widetilde{W}$. \\
(ii) Suppose that $\nabla$ is regular at $x\in M$. Then for every $w \in \widetilde{W}_{x}$
there exists a local parallel section 
$X:U\subseteq M \rightarrow \widetilde{W}$ with $X(x) = w$.
\end{theorem}
{\bf Proof:}\\
(i) follows directly from the definition of $\widetilde{W}$. \\
(ii) Suppose that $\nabla$ is regular at $x\in M$ and let $w \in \widetilde{W}_{x}$. 
By regularity, there exists a neighbourhood $U_{1}$ of $x$
and a frame $(X_{1},...,X_{n})$ of $\tilde{\pi}^{-1}(U_{1})\subseteq \widetilde{W}$.
By choosing a possibly smaller
neighbourhood $U_{2} \subseteq U_{1}$ of $x$ we can extend $(X_{1},...,X_{n})$ to a  
frame ${\cal X} := (X_{1},...,X_{n}, ...,X_{N})$ of $\pi^{-1}(U_{2})\subseteq W$. 
Let $\omega = \omega^{i}_{j}$ denote
the connection form of $\nabla$ with respect to  
${\cal X}$: $\nabla_{X}X_{j}=\sum_{i=1}^{N}X_{i} \omega^{i}_{j}(X)$.
Since $\widetilde{W}$ has zero second fundamental 1-form,
\[ \omega = \left( \begin{tabular}{c|c} 
$\phi$ & $*$ \\ \hline
$0$ & $*$  \end{tabular} \right) \]
where $\phi$ is an $n \times n$ matrix of 1-forms. The curvature
form $\Omega = \Omega^{i}_{j}$ of $\nabla$ with respect to ${\cal X}$ is

\[ \Omega = d\omega + \omega\wedge\omega = \left( \begin{tabular}{c|c}
$d\phi+\phi\wedge\phi$ & $*$ \\ \hline
$0$ & $*$ \end{tabular} \right) \]
Since the curvature tensor $R$ is identically zero, when restricted to
$TM\otimes TM\otimes \widetilde{W}$, it follows that
\[ d\phi+\phi\wedge\phi = 0\]
 
Therefore, by the Frobenius Theorem, 
there exists an $n \times n$ matrix of functions $A= A^{i}_{j}$
defined in a 
neighbourhood $U\subseteq U_{2}$ of $x$ such that 
$dA = -\phi \wedge A$ and $A(x) = I_{n \times n}$, 
the $n \times n$ identity matrix (cf. \cite{ee}, chp. 7, 2. Proposition 1., 
pg. 290). Let $c^{j}$, $1\leq j \leq n$, be real scalars satisfying 
$w = \sum_{j=1}^{n}X_{j}(x) c^{j}$. Define functions $f^{i}$ on $U$  
by 
\[ f^{i} = \left\{ \begin{array}{ll}
\sum_{j=1}^{n} A^{i}_{j}c^{j} & \hspace{0.5in} 1 \leq i\leq n \\
0 & \hspace{0.5in} n+1 \leq i \leq N. \end{array} \right. \]
Let $X:U\rightarrow \widetilde{W}$ be the local section of $\widetilde{W}$ defined  by
$X:= \sum_{i=1}^{N} X_{i}f^{i}$. Since $df+\omega\cdot f=0$, $X$ is  parallel.
Moreover, $X(x)=w$. \\
{\bf q.e.d.}

\begin{corollary} \label{corollary:flat}
Let $\nabla$ be a regular connection on the smooth vector bundle $\pi:W\rightarrow M$.
Then $(\widetilde{W}, \nabla)$ is a flat vector bundle over $M$.
\end{corollary}

\begin{corollary} \label{corollary:number}
Let $\nabla$ be a connection on the smooth vector bundle $\pi:W\rightarrow M$, regular at
$x\in M$. Then there are  $dim \hspace{0.03in} \widetilde{W}_{x}$ independent local 
parallel sections in a neighbourhood of $x\in M$. 
\end{corollary}

\hspace{-0.3in} {\bf Example} Consider the symmetric connection $\nabla$ on the 2-sphere, 
$M=S^{2}$, defined as follows:
$\Gamma^{\theta}_{\phi \phi}=-sin{\theta}cos{\theta}$, $\Gamma^{\phi}_{\theta \phi}=
\Gamma^{\phi}_{\phi \theta}=cot{\theta}$ and all other Christoffel symbols are zero. 
Here $\theta$ and
$\phi$ are the polar and azimuthal angles on $S^{2}$, respectively. Let 
\[ \begin{array}{lll}
   X_{1} & = & d\theta \otimes d\theta \\
   X_{2} & = & d\phi \otimes d\phi \\
   X_{3} & = & d\theta \otimes d\phi + d\phi \otimes d\theta 
  \end{array} \]
be a basis of $W$, the symmetric elements of $T^{*}M\otimes T^{*}M$.
The curvature terms 
$R_{\theta \phi}=\nabla_\frac{\partial}{\partial \theta}
\nabla_\frac{\partial}{\partial \phi}-\nabla_\frac{\partial}{\partial  \phi}
\nabla_\frac{\partial}{\partial \theta}$ are
\[ \begin{array}{lll}
   R_{\theta \phi}(X_{1}) & = & -(sin^{2}{\theta})X_{3}  \\
   R_{\theta \phi}(X_{2}) & = & X_{3}  \\ 
   R_{\theta \phi}(X_{3}) & = & 2X_{1}-2(sin^{2}{\theta})X_{2}  
     \end{array} \]
This gives $W^{(0)} = span(X_{1} +(sin^{2}{\theta}) X_{2})$. 
Non-zero local sections of $W^{(0)}$ are of the form 
$X=f(X_{1}+(sin^{2}{\theta}) X_{2})$ where $f$ is a smooth non-vanishing function defined 
on an open subset of 
$S^{2}$. The covariant derivative of $X$ is $\nabla X=X \otimes dlog|f|$ and so 
$W^{(1)} = W^{(0)}$. Thus $\widetilde{W} = W^{(0)}$. 
Since $\widetilde{W}$ is a rank one vector bundle over $S^{2}$  it follows that $\nabla$
is a locally metric connection; in fact, it is the Levi-Civita connection of the induced
metric of the standard embedding of the two-sphere in three-dimensional Euclidean space.

\end{document}